\documentclass[12pt]{amsart}
\usepackage{amscd,amsthm,amsmath} 
\newtheorem{thm}{Theorem}[section]
\newtheorem{cor}[thm]{Corollary}
\newtheorem{lem}[thm]{Lemma}
\newtheorem{prop}[thm]{Proposition}

\theoremstyle{definition}
\newtheorem{defn}[thm]{Definition}
\theoremstyle{remark}
\newtheorem{rem}[thm]{Remark}
\numberwithin{equation}{section}

\newcommand{\A}{\mathcal{A}}

\newcommand{\N}{{\mathbb N}}
\newcommand{\Z}{{\mathbb Z}}
\newcommand{\Q}{{\mathbb Q}}

\begin{document}

\title {Secondary Brown-Kervaire Quadratic forms and $\pi$-manifolds}
\author{Fuquan Fang}
\address{ Nankai Institute of Mathematics, Nankai University,
Tianjin 300071, P.R.C }

\email{ ffang@sun.nankai.edu.cn}
\author{Jianzhong Pan }%
\address{Institute of Math.,Academia Sinica ,Beijing 100080 ,China
and \newline Department of Mathematics Education , Korea University , Seoul , Korea }%
\email{pjz62@hotmail.com}%

\thanks{The first author was supported in part by NSFC 1974002, Qiu-Shi Foundation and CNPq
and the second author is partially supported by the
NSFC project 19701032 and ZD9603 of Chinese Academy of Science and
Brain Pool program of KOSEF }%

\date{Dec. 20,1999}


\begin{abstract}
In this paper we assert that for each $\Phi$-oriented
$2n$-manifold (c.f : Definition 1.1) $M$ where $n\ge 4$ and
$n\ne 3(mod 4)$, there is a well-defined quadratic function
$\phi_M: H^{n-1}(M, \Z_4)\to \Q/\Z$, we call the secondary Brown-Kervaire quadratic forms, so that
\begin{itemize}
\item{ $\phi _{M}(x+y)=\phi _{M}(x)+\phi _{M}(y)+j(x\cup Sq^2y)[M]$},
\item{ the Witt class of $\phi _M$ is a homotopy invariant, if
the Wu class $ v_{n+2-2^i}(\nu _M)=0$ for all $i$.}
\end{itemize}
 where $j: \Z_2 \to \Q/\Z$ is the inclusion homomorphism and
$\nu _M$ the stable normal bundle of $M$.

Among the applications we obtain a complete classification of
$(n-2)$-connected $2n$-dimensional $\pi$-manifolds up to
homeomorphism and homotopy equivalence, where $n\geq 4$ and
$n+2\neq 2^i$ for any $i$. In particular, we prove that the
homotopy type of such manifolds determine their homeomorphism
type.

\end{abstract}
\maketitle

\section{Introduction}\label{S:intro}

Let $M$ be a $2n$-dimenisonal framed manifold (i.e. a $\pi$-manifold with a framing) where $n=1(mod2)$.
The Kervaire invariant of $M$ is the Arf invariant of a $\Z_2$-valued Kervaire quadratic form of $M$
$$q_M: H^{n}(M, \Z_2) \to  \Z_2$$
satisfying
$$q_M(x+y)=q_M(x)+q_M(y)+(x\cup y)[M]_2 \hspace{2cm} (1.1) $$
It was invented by Kervaire to find the first example of non-smoothable PL-manifold. Kervaire invariants
and its various generalizations, e.g. the Brown-Kervaire invariants[4], play very important roles in
geometric topology. Formally, $q_M$ is a ``quadratic form'' subject to the symmetric bilinear form
$$
\begin{array}{cc}
H^n(M,\Z_2)\times H^n(M,\Z_2)\to \Z_2\\
(x,y)\to x\cup y[M]_2
\end{array}
$$

For a Spin manifold of even dimension, there is another symmetric bilinear form $\mu _M$ studied by
Landweber and Stong [16]:
$$
\begin{array}{cc}
\mu _M: H^{n-1}(M,\Z_2)\times H^{n-1}(M,\Z_2)\to \Z_2\\
(x,y)\to Sq^2(x)\cup y[M]_2
\end{array}
$$
A natural algebraic question to ask is whether there is an intrinsic ``quadratic form'' of $M$ subject to
$\mu _M$.
To answer this turns out to be the main novelty of this paper. For a large family of Spin manifolds
including all $\pi$-manifolds, the so called
$\Phi$-oriented manifolds, we will define a $\Q/\Z$-form subject to $\mu _M$, which resembles to the
Brown-Kervaire quadratic forms in the formulation. It has the most similar properties of the Brown-Kervaire
quadratic forms, e.g., the isomorphism class of the form is a homotopy invariant if the manifold has
vanishing Wu classes. A bit surprising to us, this invariant applies to give a classification of
$(n-2)$-connected $2n$-dimensional $\pi$-manifolds up to homotopy equivalence and homeomorphism ($n\ge 4$).

To state our main results, let us start with some notations.

Let $\{ Y_k\}_{k\in  \N} $ be a connected spectrum with $U\in H^{0}(Y)\cong  \Z$ a generator so that
$i^{*}U\in H^0(S^0)$ a generator, where $i: S^0 \to Y$ is the inclusion map of the spectrum.

\begin{defn}
 (i) $\{ Y_k\}_{k\in  \N} $ is called $\Phi$-{\it
orientable} if $Sq^2U=0$, $\chi (Sq^{n+2})(U)=0$ and $0\in \Phi
(U)$, where $\Phi$ is a secondary cohomology operator associated
with the Adem relation (see Section \ref{S:2} for the definition):
$$
\begin{array}{ll}
\chi(Sq^n)Sq^3+\chi(Sq^{n+2})Sq^1+Sq^1\chi(Sq^{n+2})=0   & n=2(mod
4) \\ \chi(Sq^n)Sq^3+Sq^1\chi(Sq^{n+2})=0   &n=0(mod 4) \\
\chi(Sq^{n+1})Sq^2+Sq^1\chi(Sq^{n+2})=0  & n=1(mod 4)
\end{array}
$$ where $\chi: \A_2\to \A_2$ is the anti-automorphism of the
Steenrod algebra $\A_2$ \cite{adams}.
\end{defn}

A spherical fibration $\xi$ (a manifold) is called  $\Phi
$-{\it orientable} if its Thom spectrum $T\xi $ (stable normal bundle
$\nu _M$) is. We define the {\it universal} $\Phi $-orientable $\Omega$-spectrum
$\widetilde W(n)$ by setting  $\widetilde {W}_k{(n)}$ to be the total space of the following
Postnikov tower:
\newline
$$
\begin{array}{ccccc}
&  & \widetilde{W}_{k}(n) &  &  \\ &  & \downarrow {\Pi _2} &  &
\\ &  & W_{k}(n) & \stackrel {k_2}{\longrightarrow} & K_{k+n+2}
\\ &  & \downarrow   {\Pi _1} &  &  \\ &  & K(\Z,k) &
\stackrel {Sq^2\times \chi (Sq^{n+2}) } {\longrightarrow} &
K_{k+2}\times K_{k+n+2}
\end{array}
$$
where $K_i=K(\Z_2, i)$, $K(\Z, i)$ are the Eilenberg-Maclane spaces, $k_2\in \Phi ({\Pi _1 } ^*l_{k})$
and $l_k$ is the basic class.

Note that a spectrum $Y$ is $\Phi$-orientable if and only if $U\in H^0(Y)$ can be lifted to a map
$w: Y\to \widetilde W(n)$. We call such a lifting a $\Phi$-{\it orientation} of $Y$. A $\Phi$-orientation of
a manifold is understood as a $\Phi$-orientation of its Thom spectrum.

\begin{rem}\label{T:1.2}
The sphere spectrum $S^0$ is $\Phi $-orientable. Thus
stably parallelizable manifolds are $\Phi$-orientable.
\end{rem}

Our main results are:

\begin{thm}\label{T:0.6}
{\it Let $M$ be a $\Phi$-oriented manifold of dimension $2n$, where $n\ne 3(\mbox{mod } 4)$. Then there is
a  function $\phi_{M}: H^{n-1}(M, \Z_4)\to \Q/\Z$ such that, for all $x, y\in H^{n-1}(M,\Z_4)$,
 $$\phi _{M}(x+y)=\phi _{M}(x)+\phi _{M}(y)+j(x\cup Sq^{2}y)[M],$$ where
$j:  \Z_2 \to \Q/\Z $ is the inclusion.}
\end{thm}

\begin{rem} In general, $\phi _M$ depends on the $\Phi$-orientation, just like the Kervaire quadratic
form depends on the framing of the manifold. We will prove that $\phi _{M}(x)$ depends only
on the $\Phi $-oriented bordism class $[M,x]$.
\end{rem}

\begin{rem} If $n=3(\mbox{mod }4)$, the analogous definition gives only a linear function.
\end{rem}

Let $BSpin _G$ be the classifying space for spherical Spin fibrations. By Brown[4], a {\it Wu orientation }
of a Spin spherical fibration $\xi \searrow M$ is a lifting of the classifying map $\xi : M\to BSpin_G$ to
$BSpin_G\langle v_{n+2}\rangle $. A {\it Wu orientation } of $\nu _M$, the
stable normal bundle of $M$, is understood as a Wu orientation of
$M$, where $BSpin_G\langle v_{n+2}\rangle \to BSpin_G$ is a
principal fibration with $v_{n+2}\in H^{n+2}(BSpin_G, \Z_2)$ as the
$k$-invariant.

We call quadratic forms $\phi _{M_i}: H^{n-1}(M_i, \Z_4)\to
{\Q/\Z}$, $i=1,2$  {\it Witt equiavlent} if there exists an
isomorphism $\tau : H^{n-1}(M_1, \Z_4)\to H^{n-1}(M_2, \Z_4)$ so
that $\phi _{M_2}(\tau x)=\phi _{M_1}(x)   $ for all $x\in
H^{n-1}(M_1, \Z_4)$.

\begin{thm}\label{T:0.11}
 {\it Let $M_1$ and $M_2$ be $\Phi$-oriented $2n$-manifolds.
Suppose that the Wu classes $v_{n+2-2^j}(\nu _{M_i})=0$ for all $2^j\le n+2$. If $f: M_1\to
M_2$ is a homotopy equivalence preserving the spin structure
(resp. Wu orientation) if $n=0, 1(mod 4)$ (resp. $n=2(mod 4)$).
Then $$\phi _{M_1}(f^*x)=\phi _{M_2}(x)$$ for all $x\in
H^{n-1}(M_2, \Z_4)$.}
\end{thm}

Since the Wu class $v_0=1$, the assumption in the above theorem implies that $n+2\ne 2^i$ for any integer
$i$.

For framed manifolds, the Brown-Kervaire secondary quadratic forms have the following property:

\begin{prop}\label{T:0.9}
If $M$ is a framed manifold of dimension $2n$, where $n\ne 3(mod 4)$. Then $\phi _M$
factors through $\Z_4\subset {\Q/\Z}$ $(\text{resp. } \Z_2 \subset {\Q/\Z})$,
provided $n=2(mod4)$ $(\text{resp. }n=0, 1(mod4))$.
\end{prop}

To state the next results, we need some preliminaries.

Let $H$ be a finitely generated abelian group, and
 $$\begin{array}{rll}\mu : Hom(H, \Z_2 ) \otimes Hom(H,  \Z_2)
 \to \Z_2
\end{array}$$
be a symmetric bilinear form.  We say that $\mu $ is of \underline
{\it diagonal zero} if $\mu (x,x)=0$ for each $x \in Hom(H, \Z_2)$.
A function $\phi : Hom(H,  \Z_4 ) \to  \Q/\Z$ is called \underline
{\it quadratic} with respect to $\mu $ if $$
\begin{array}{rll}
\phi (x+y) = \phi (x) + \phi (y) + j( \mu (x, y))
\end{array}$$
where $j :  \Z_2 \to \Q/\Z $ is the inclusion. This gives a triple
$(H, \mu , \phi )$.  We say triples $(H_1, \mu _1, \phi _1)$, and
$(H_2, \mu _2, \phi _2)$ are {\it isometric} if there exists an
isomorphism $\tau: H_1\to H_2$ such that $\mu _1(x, y)=\mu _2(\tau
x, \tau y)$ and $\phi _1(x)=\phi _2(\tau x)$ for all $x, y$. We
denote by $[H, \mu , \phi ]$  the isometry class of  a triple.

\begin{rem}
Since the natural map $ Hom(H_{n-1}(M),  \Z_2 ) \to H^{n-1}(M,
\Z_2)$ is not an isomorphism in general, the notions of isometry
associated with $\mu$ and $\mu_M$ as above are different. They do
agree however for $(n-2)$-connected manifolds which we will assume
in the later application. We will use both of them when necessary.
\end{rem}

Let ${i}$ denote the maximal exponent of the $2$-torsion subgroup of
$H_{n-1}(M)$ and let $Sq_{i}^1 \in H^n(K( \Z_{2^{i}}, n-1),
\Z_2)\cong \Z_2$ be the unique generator. Considering $Sq_{i}^1 $ as
a cohomology operation we get a function $$\begin{array}{rll}
q_{M}(Sq_{i}^1): H^{n-1}(M, \Z_{2^i}) \to \Z_2.
\end{array}$$
This gives a homomorphism since $Sq^1_ix\cup Sq^1_iy=Sq_i^1(x\cup
Sq_i^1y)=0$ for $x, y\in H^{n-1}(M, \Z_2)$.   We denote by
$[H_{n-1}(M), \mu _M, q_M(Sq^1_i)]$ for the isometry class of the
triple.  By \cite{browd2}, the Kervaire invariant of a smooth
framed manifold of dimension $2n$, where $n\ne 2^i-1$, is zero.
For $i\leq 5$, there are smooth manifolds of dimension $2^{i+1}-2$
of Kervaire invariant $1$. It is still an open problem whether
there is such a manifold for  $i\geq 6$.

The Kervaire invariant does not depend on the framings
of the underlying $2n$-manifold if $n\ne 1, 3, 7$ and the manifold
is highly connected, e.g. $(n-2)$-connected. Moreover, by
\cite{brown} the Kervaire form is a homotopy invariant if $n\ne 1, 3, 7$ and $(n-2)$-connected.

Let $M$ be a $(n-2)$-connected $2n$-dimensional $\pi$-manifold.
Observe that if $n\ge 3$, there exists a $(n-2)$-connected
$\pi$-manifold, $N$, so that $M=N\# X$ and $H_{n}(N, \Q)=0$, where
$X$ is a $(n-1)$-connected $2n$-manifold. Since the classification
of $(n-1)$-connected $2n$-manifolds is well understood
\cite{wall2}, for convenience in the following theorem we
assume that $H_{n}(M, \Q)=0$. For such a manifold, consider the
correspondence $$\pi: M \longmapsto [H_{n-1}(M), \mu _{M}, \phi
_M] (\mbox{ resp.}\text{ } [H_{n-1}(M), \mu _{M}, \phi _M,
q_M(Sq^1_{2^i})])$$ if $n=0(mod 2)$ (resp. $n=1(mod 2)$).

In the following theorem let $\alpha(n+2)$ be the number of $1's$
in the binary expansion of $n+2$.

\begin{thm}\label{T:0.13}
{\it Suppose $n\geq 4$ and $\alpha (n+2)\geq 2$. Then $\pi$ gives
a 1-1 correspodence between the homeomorphism types (resp.
homotopy types) of $(n-2)$-connected $2n$-dimensional
$\pi$-manifolds $M$ so that $H_n(M, \Q)=0$ with the following
algebraic data

\noindent  (a)  $\wp _n=\{ [H, \mu , \phi ]: \mbox{ diag }\mu =0
\mbox{ and } \phi \mbox{ factors through } j: \Z_4\to \Q/\Z\}$
\mbox{ if $n=2(mod 4)$},\newline (b) $\wp _n=\{ [H, \mu , \phi ]:
\phi \mbox{ factors through }j: \Z_2\to \Q/\Z\}$ if $n=0(mod
4)$,\newline (c)  $\wp _n=\{ [H, \mu , \phi , \omega ]: \omega \in
Hom(tor (H)\otimes \Z_{2^i}, \Z_2),  \phi $ factors through $j:
\Z_2\to \Q/\Z \}$ if $n=1(mod 4)$,\newline (d) $\wp _n=\{ [H, \mu
, \omega ]:  \omega \in Hom(tor H\otimes
 \Z_{2^i}, \Z_2)\}$ if $n=3(mod 4)$.

\noindent  where $i$ is the highest exponent of the $2$-cyclic
subgroup of $H$ and if $n=1(mod 2)$, the pairing $\mu (x, x)=0$
(resp. $\delta \omega (x)$) if $x$ can be lifted to a $\Z_4$ class
with order $4$ (resp. $x$ is of order $2$), $\delta \in \{ 0, 1\}$
is ambiguous. }
\end{thm}

\begin{rem}
The classification of $(n-2)$-connected $2n$-manifolds with torsion
free homology groups has been given by Ishimoto[9][10]. But his method does not work
if the homology group has torsion.
\end{rem}

The organization of this paper is as follows.

In $\S$\ref{S:1} we define the secondary Brown-Kervaire form and state its basic properties.

In $\S$\ref{S:2}, we set up the necessary foundations on the stable homotopy theory of the
Eilenberg-Maclane spaces.

In $\S$\ref{S:3}, we are addressed to show Theorems 1.3 and 1.6.

In $\S$5, we prove Theorem \ref{T:0.13}.

\bigskip

\section{A $\Q/\Z$-quadratic form of $\Phi$-oriented manifolds}\label{S:1}

Let us begin with some conventions.  All homology/cohomology
groups will be with integral coefficients unless otherwise stated.
All spaces will have base points. Let\newline (i) $[X, Y]$ denote
the set of homotopy classes of pointed maps from $X$ to $Y$.
\newline (ii) $\{ X, Y\} =lim[S^kX, S^kY]$.\newline (iii)
$\pi_*^s(X)$  be its 2-localization to simplify the notation.
\bigskip

Let $\kappa : K(\Z_4, n-1) \times  K(\Z_4, n-1) \to  K(\Z_4, n-1)$ be
the multiplication of  $K(\Z_4, n-1)$ and let $H(\kappa )$ be the
Hopf construction of $\kappa$.

\begin{prop}\label{T:0.2}
{\it The homomorphism $$H(\kappa )_*: \pi ^s_{2n}(K(\Z_4, n-1)\land
K(\Z_4, n-1))\to \pi_{2n}^s(K(\Z_4, n-1))  $$ is injective  if
$n\neq 3(mod 4)$, and zero if $n=3(mod4)$.}
\end{prop}
\begin{rem}
If $\Z_4$ is replaced by $\Z_2$, then $H(\kappa )_*$
is trivial.
\end{rem}

By Theorem \ref{T:stable} and the proof of it, we obtain
$$\begin{array}{lcc} \pi ^s_{2n}(K(\Z_4, n-1)\land K(\Z_4,
n-1))\cong \Z_2\mbox{ if $n\geq 4$,}\\ \pi ^s_{2n}(K(\Z_4,
n-1))\cong \Z_4 \mbox{ if $n=2(mod 4)$.}
\end{array}
$$

Let $\lambda _0$ be a generator of Im $H(\kappa )_*$ if $n\neq
2(mod4)$, and a specified generator of $\pi ^s_{2n}(K(\Z_4,
n-1))\cong \Z_4$ otherwise. For a given spectrum $Y$, let
$$H_{*}(K(\Z_{4}, n-1); Y)=\mbox{lim }\pi _{*+k}
(K(\Z _4, n-1)\wedge Y_k).$$

\begin{thm}\label{T:0.3}
 {\it Suppose that $\{Y_k\}_{k\in \N}$ is a $\Phi$-orientable
spectrum. Then there exists a homomorphism $$h:   H_{2n}(K(\Z_4,
n-1); Y) \to  \Q/\Z $$ such that $h(\lambda)=\frac {1}{4}$
$(\text{resp. } \frac{1}{2})$ if $n=2(mod 4)$ $(\text{resp. }n=0,
1(mod 4))$, where $\lambda =i_*(\lambda_0)$ and $i_*: H_{2n}(K(\Z_{4},
n-1); S^0)\to H_{2n}(K(\Z_{4}, n-1); Y) $ is induced by the inclusion.}
\end{thm}

\begin{defn} A Poincar\'e triple $(M, \xi , \alpha )$ of
dimension $2n$ consists of \newline (i) A CW complex $M$ with
finitely generated homology.\newline (ii) A fibration $\xi $ over
$M$ with fiber homotopy equivalent to $S^{k-1}$, $k$ large.
\newline (iii) $\alpha \in \pi _{2n+k}(T\xi)$ such that an $(2n+k)$
Spanier-Whitehead S-duality is given by
$$S^{2n+k}  \stackrel {\alpha }{\longrightarrow} T\xi
\stackrel {\Delta }{\longrightarrow}T\xi\land M^+$$
where $T\xi$ is the Thom complex of $\xi$ and $\Delta $ is the diagonal map.
\end{defn}

Let $A_{\alpha}: \{M_{+}, K(\Z_4, n-1)\} \to \{S^{2n+k}, T\xi\wedge
K(\Z_4, n-1)\}$ be the $S$-duality map.

\begin{defn}
Let $(M, \xi, \alpha)$ be a Poincar\'e triple and
$w$ is a $\Phi$-orientation of the Thom spectrum $T\xi $. For a homomorphism $h$ in Theorem 2.3,
let $$\phi_{w,h} : H^{n-1}(M,  \Z_4)\to \Q/\Z$$ be defined by setting
$$\phi
_{w,h}(x)=h([(w\wedge id) A_{\alpha}(x)]).
$$
\end{defn}

\begin{thm}
{\it Let $\phi _{w,h}$ be defined as above. Then for all $x, y\in H^{n-1}(M,\Z_4)$,  \newline
(i) If $n\ne 3(\mbox{mod } 4)$, the function is quadratic, i.e.
$$\phi _{w,h}(x+y)=\phi _{w,h}(x)+\phi _{w,h}(y)+j(x\cup Sq^{2}y)[M]$$
where $j:  \Z_2 \to \Q/\Z $ is the inclusion;\newline
(ii) If $n=3(mod 4)$, $\phi _{w,h}$ is linear, i.e.
$$\phi _{w,h}(x+y)=\phi _{w,h}(x)+\phi _{w,h}(y).$$}
\end{thm}

Now we want to study how the function $\phi_{w,h}$ depends on the
choice of  the orientation of the Thom spectrum $T\xi$.

Let $w_i$, $i=1, 2$ are orientations of the Thom spectrum $T\xi$.
Let $$d_{1}(w_1, w_2)\in H^1(T\xi)\oplus H^{n+1}(T\xi )$$ denote
the difference of the composition maps $\Pi _2$$w_1$ and
$\Pi_2$$w_2$, where $\Pi _2$ is as in the definition of the
universal $\Omega$-spectrum $\widetilde W(n)$. Clearly, $w_1$ and
$w_2$ are homotopy if and only if $d_1(w_1, w_2)=0$ and a
secondary obstruction vanishes. The following theorem shows that
the secondary  obstruction does not affect our quadratic function
$\phi _{w,h}$.

\begin{thm} Let $\phi _{w_i,h}$ be the quadratic forms associated with $(w_i, h)$, $i=1, 2$.
If $d_1(w_1, w_2)=0$, then $\phi _{w_1,h} (x) =  \phi _{w_2,h}(x) $ for all
$x\in H^{n-1}(M, \Z_4)$.
\end{thm}

In general, the quadratic form $\phi _{w,h}$ does depend on the choice of $w$ and $h$. In order to obtain a
well-defined invariant of the $\Phi$-oriented manifold, we now choose certain type of $\Phi$-orientations
of the Thom spectrum $T\xi$ in an universal way and then define the Brown-Kervaire secondary quadratic forms
to be the quadratic functions associated to those $\Phi$-orientations.

Let $\gamma  \searrow BSpin_G $ be the universal Spin spherical fibration and
$U\in H^0(MSpin_G, \Z_2)$ the universal Thom class. Note that
$$\begin{array}{lccc} \chi
(Sq^{n+2})U =\chi (Sq^{n+1})Sq^1U =0\mbox{ if $n$ is odd}\\ \chi
(Sq^{n+2})U =\chi (Sq^n)Sq^2U =0\mbox{ if $n=0(mod 4)$},
\end{array}
$$
Thus $U$ lifts to a map $f: MSpin_G\to W (n)$. By the Thom isomorphism, $f^*k_2$ gives an element of
$\bar {k_2} \in H^{n+2}(BSpin_G, \Z_2)$. Let $\pi :
BSpin_G\langle \bar{k_2}\rangle \to BSpin_G $ be the principal fibration with $k$-invariant
$\bar{k_2}$.

If $n=2(mod 4)$, we get a similar principal fibration $\pi :
BSpin_G\langle \bar{k_2}\rangle \to BSpin_G \langle v_{n+2}\rangle$,  where $BSpin_G\langle
v_{n+2}\rangle \to BSpin_G $ is the fibration with fibre $K_{n+1}$
and $k$-invariant $v_{n+2}$.

It is easy to see that the fibration $\pi ^*\gamma $ is $\Phi$-orientable. Clearly  the classifying map
of every $\Phi$-orientable stable spherical fibration lifts
to $BSpin_G\langle \bar{k_2}\rangle $.

\begin{defn}
The fibration  $\pi ^*\gamma $ is called the {\it universal $\Phi$-orientable spherical Spin
fibration}. Its Thom spectrum,   $MSpin_G\langle \bar{k_2}\rangle$, is called the  {\it universal
$\Phi$-orientable Thom spectrum}.
\end{defn}

For a closed $\Phi$-orientable manifold $M^{2n}$, there is a Poincar\'e triple $(M, \nu _M, \alpha )$
where $\nu _M$ is the stable normal bundle and $\alpha \in \pi _{2n+k}(T\nu _M )$ is the normal invariant
of $M$ (obtained by the Thom-Pontryagin construction.)

\begin{defn}
Fix a  connected spectral map ${\bf u}:  MSpin_G\langle
\bar{k_2}\rangle \to \widetilde {W}(n)$ and a homomorphism $h$ in
Theorem 2.3. For a $\Phi$-orientable manifold $M$, let $$\phi
_M=\phi _{w,h}$$ where $w={\bf u}\circ T(v)$ and $T(v)$ is the
Thom map of a classifying bundle map of the stable bundle $\nu
_M$.
\end{defn}

Now we prove Theorem 1.6 assuming Theorem 2.7.

\begin{proof}[Proof of Theorem \ref{T:0.11}]
Let $\xi _i=\nu _{M_i}$ be the stable normal bundle of $M_i$ and $\alpha _i\in \pi _{2n+k}(T\xi _i)$ be the
normal invariant, $i=1,2$. By the definition, $\phi _{M_i}=\phi _{w_i, h}$ where $w_i={\bf u}\circ T(v_i)$
and $T(v_i): T(\xi _i)\to MSpin _G\langle \bar k_2\rangle$ the Thom map.

Let $\tilde f: f^*\xi _2\to \xi _2$ be a bundle map over the homotopy equivalence
$f$. Let $\alpha _3=T(\tilde f)^{-1}_*\alpha _2$, where $T(\tilde f)$ is the Thom map of $\tilde f$.
The Poincar\'e triple $(M_1, f^*\xi _2, \alpha _3)$ together with the $\Phi$-orientation
$w_2\circ T(\tilde f)$ gives a quadratic form $\phi _3$, where $w_2={\bf u}\circ T(v_2)$ is a
$\Phi$-orientation of $M_2$. By 2.5 we get that
 $$\phi _3(f^*x)=\phi _{M_2}(x)$$
 for all $x\in H^{n-1}(M_2, \Z_4)$.

To prove the desired result, it suffices to prove $\phi _3=\phi _{M_1}$.

Note that $f^*\xi _2$ and $\xi _1$ are stably equivalent as
spherical fibration since $f$ is a homotopy equivalence. Thus we can regard $f^*\xi _2$ and $\xi _1$ as the
the same and so get two orientations for $\xi _1$, $({\bf u}\circ T(v_1), h)$ and $({\bf u}\circ T(v_2)
\circ T(\tilde f), h)$.  Since $f$ preserves the Spin structures/Wu orientations, $\pi \circ
v_2\circ f\simeq \pi \circ v_1$, where $\pi: BSpin_G \langle \bar
{k_2}\rangle \to BSpin_G$/$BSpin_G \langle v_{n+2} \rangle$ is the
principal fibration as above. This clearly implies that there exists a fibre automorphism
$g\in Aut (\xi _1)$ over the identity such that
$$T(\pi \circ v_2\circ \tilde f) \simeq
T(\pi \circ v_1)\circ T(g).$$
 Notice that $g$ gives a unique element $g _0\in [M_1, G_k]$, where $G_k$ is the space of self homotopy
equivalences of $S^k$.  By a formula in Brown \cite{brown}, the $(n+1)$-dimensional component
of $d_1({\bf u}\circ T(v_1)\circ T(g),{\bf u}\circ T(v_1))$ is
$\sum v_{n+2-2^i}\cup g_0^*u_{2^i-1}$, where $u_{2^i-1}$ is the transgression of $w_{2^i}\in
H^{2i}(BG_k, Z_2)$.  By assumption, it must vanish since the Wu classes vanish.
On the other hand, the $1$-dimensional component of  $d_1({\bf u}\circ T(v_1)\circ T(g),{\bf u}\circ T(v_1))$
is determined by the Spin structures and so it vanishes since $f$ preserves the Spin structures.
By Theorem 2.7 it follows that
$$\phi _{M_1}=\phi _4,$$
the quadratic form associated with the Poincar\'e triple  $(M_1, \xi _1, \alpha _1)$ and the
$\Phi$-orientation $w_2\circ T(\tilde f)$.

 Note that in the definitions of $\phi _3$ and $\phi _4$
the only different ingradients are the normal invariants, after identifying $\xi _1$ with
$f^*\xi _2$. By Theorem 2.7 once again   $\phi _3=\phi _4$. This implies the desired result.
\end{proof}

Now we prove Proposition 1.7.

\begin{proof}[Proof of Proposition 1.7]
Since $M$ is a framed manifold, the stable normal bundle is trivial, i.e. the classifying map
of $\nu _M$ factors through a point. Choose a $\Phi$-orientation $w={\bf u}\circ T(v): \nu _M$ with
$v$ the bundle map of $\nu _M$ to the trivial $k$-bundle on a point, then $\phi _M(x)$ factors through
the stable homotopy group $\pi _{2n}^s(K(\Z_4, n-1))$. By Theorem 3.1  $\pi _{2n}^s(K(\Z_4, n-1))\cong \Z_4$
if $n=2(mod 4)$ and the order of elements in  $\pi _{2n}^s(K(\Z_4, n-1)$ is at most $2$ if
$n=0, 1(mod 4)$. On the other hand, by Theorem 1.6 the definition of $\phi _M$ does not depend on the choice of the $\Phi$-orientations since $M$ is a framed manifold. This completes the proof.
\end{proof}

\bigskip

\section{Some preliminaries on stable homotopy theory}\label{S:2}

In this section we calculate  the stable homotopy groups $\pi
^s_{2n}(K(\pi , n-1))$ (see Theorem \ref{T:stable}).  We will also
introduce some  $2$-stage Postnikov tower which will give  the
secondary cohomology operation $\Phi$ used in Section
$\S$\ref{S:intro}.

\begin{thm}\label{T:stable}
 {\it The $2n$-th stable homotopy group of $K(\pi , {n-1})$
for $n\geq 4$ is as follows:

{\small \begin{tabular}{|r|c|c|c|c|c|}\hline
$n\geq 4$  & $0(mod 4)$ & $1(mod4)$  & $2(mod4)$ & $3(mod4)$
\\ \hline
$\pi_{2n}^s(K(\pi ,n-1))$  & $ (\Z_2)^{2(t+k)+s+p}   $ & $
(\Z_2)^{t+2k+s+p}   $ & $ (\Z_4)^{t+k}\oplus (\Z_2)^{s+p}   $& $
(\Z_2)^{k+s+p}   $ \\ \hline
\end{tabular} }
\newline where $p={{t+k+s}\choose{2}}$ and $\pi =G_0\times \Z^t\times
\Z_{2^{i_1}} \times \cdots \times \Z_{2^{i_k}} \times \Z_2^s$,
$i_j\geq 2$ if $1\leq j\leq k$ and $G_0\otimes \Z_2=0$.}

When $\pi=\Z$, Theorem \ref{T:stable} follows from
\cite{mahowald1}.
\end{thm}

\begin{proof}
It is easy to know that(since we are computing the 2-localization)
\[
\pi ^s_{2n}(K(\pi , n-1)) =\pi ^s_{2n}(K(\pi/G_0 , n-1))
\]
Assume $G_0=0$ from now on. If $\pi=\pi_1 \bigoplus \pi_2$ with
$\pi_1$ nontrivial and $\pi_2$ a nontrivial cyclic group, then
\[
K(\pi , n-1)=K(\pi_1 , n-1) \times K(\pi_2 , n-1)
\]
and we have by a result in \cite{barcus} that
\[
\pi ^s_{2n}(K(\pi , n-1))=\bigoplus_{i=1,2}\pi ^s_{2n}(K(\pi_i ,
n-1)) \bigoplus \pi ^s_{2n}(K(\pi_1 , n-1)\wedge K(\pi_2 , n-1))
\]
\[
=\bigoplus_{i=1,2}\pi ^s_{2n}(K(\pi_i , n-1)) \bigoplus
H_{n+1}(K(\pi_1 , n-1),\pi_2)
\]
An easy calculation shows that $ H_{n+1}(K(G_1 , n-1),G_2)=Z_2$ if
$G_1,G_2$ are nontrivial cyclic groups and thus $ H_{n+1}(K(\pi_1
, n-1),\pi_2)=Z_2^{t+k+s-1}$.

On the other hand we know groups $\pi ^s_{2n}(K(\Z , n-1))$ and
$\pi ^s_{2n}(K(\Z_2 , n-1))$ by results in
\cite{milg},\cite{mahowald1}. To complete the proof it remains to
calculate $\pi ^s_{2n}(K(\Z_{2^i} , n-1))$ for $2\leq i < \infty$
which will be given in the following results.
\end{proof}

Recall that for each locally finite connected CW complex $X$ we
can define a space $$\Gamma _qX= S^{q-1}\propto_{T}X\land X=
S^{q-1}\times (X\land X) /\{ (x,y,z)\sim (-x,z,y); (x, *)\sim *
\}$$ for every $q\in Z_+$. By \cite{milg} Theorem 1.11, for a
$(n-2)$-connected space $X$, $\Gamma _qX$ is $(2n-3)$-connected.
Moreover, if $X=K(\pi , n-1)$,  we have a fibration $$G_q\to
\Sigma ^qK(\pi , n-1) \to K(\pi , q+n-1)$$ where $G_q\simeq \Sigma
^q\Gamma _q(K(\pi , n-1))$ through dimension $(3n+q-3)$.  Thus
$\pi _i^s(K(\pi , n-1))\cong  \pi _i^s(\Gamma _q(K(\pi , n-1))$
for $n< i<3n-3$.

When $q=1$, $\Gamma _qX= X\land X$ .The corresponding sequence is
:
\[
\Sigma F_{n-1}(\pi) \overset{H(\kappa)}{\to} \Sigma K(\pi , n-1)
\to K(\pi , n)
\]
where $F_{n-1}(\pi)=K(\pi , n-1) \land K(\pi , n-1)$. After
 $q-1$ time suspensions we get a fibration sequence at least in
dimensions less than $3n+q-4$
\[
\Sigma^q F_{n-1}(\pi) \overset{\Sigma^{q-1}H(\kappa)}{\to}
\Sigma^q K(\pi , n-1) \to \Sigma^{q-1}K(\pi , n)
\]
Let $q$ be large enough so that we are always in the stable range
and let $r=q+2n$ ,then we have an exact sequence
\[ \label{T:sequence}
\cdots \to \pi ^s_{2n+2}(K(\pi , n))\overset{\partial}{\to}\pi
_{r}(\Sigma^qF_{n-1}(\pi))\to
\]
\[
 \to \pi ^s_{2n}(K(\pi ,n-1))\to
\pi ^s_{2n+1}(K(\pi , n))\overset{\partial}{\to} \cdots
\]
Since we know that $\pi _{r}(\Sigma^qF_{n-1}(\pi))=\Z_2$ for
$\pi=\Z$ and $\Z_{2^i}$ , we can determine inductively $ \pi
^s_{2n}(K(\pi ,n-1))$ up to extension if we know the map
$\partial$.
\begin{lem}\label{T:2.2}
For $\pi=\Z$ or $\Z_{2^i}$, a homotopy class $[g]\in \pi
^s_{2n+2}(K(\pi ,n))$ has $\partial g \neq 0$ iff
$g^*(\Sigma^{q-1}(\iota \cup Sq^2\iota))\neq 0$ where $g:S^{r+1}
\to \Sigma^{q-1}K(\pi,n)$,$\iota \in H^n(K(\pi,n),\Z_2)$ is the
generator and $g^*:H^*(\Sigma^{q-1}K(\pi,n),\Z_2) \to
H^*(S^{r+1},\Z_2)$
\end{lem}
The proof is similar to that of Lemma 1.3 in \cite{mahowald1}. The
key points are the followings:
\begin{itemize}
\item{$g^*$ can be nonzero only on element $\Sigma^{q-1}(\iota \cup Sq^2\iota)$}
\item{the Hurewicz homomorphism $H:\Z_2 \cong\pi
_r(\Sigma^qF_{n-1}(\pi)) \to H_r(\Sigma^qF_{n-1}(\pi)) $ is
nonzero}
\end{itemize}
The first statement is clear while the second is an easy
consequence of the Whitehead exact sequence(c.f. \cite{white},
page 555)

With the lemma above we can now prove the Proposition\ref{T:0.2}.
\begin{proof}[Proof of Proposition\ref{T:0.2}]
It suffices to prove that $\partial$ is trivial if $n=
0,1,2(mod4)$ and nontrivial if $n= 3(mod4)$.

For $n= 0,1,2(mod4)$ there is no $g:S^{r+1} \to
\Sigma^{q-1}K(\Z,n)$ such that $g^*(\Sigma^{q-1}(\iota \cup
Sq^2\iota))\neq 0$ since $\Sigma^{q-1}(\iota \cup Sq^2\iota)$ is
detected by the secondary cohomology operation $\varphi_n$ in
\cite{mahowald2}. Thus there is no $g:S^{r+1} \to
\Sigma^{q-1}K(Z_{2^i},n)$ such that $g^*(\Sigma^{q-1}(\iota \cup
Sq^2\iota))\neq 0$ by the naturality of  secondary cohomology
operation and the fact that $\rho_{2^i}:K(\Z,n)\to K(\Z_{2^i},n)$
corresponding to mod $2^i$ reduction induces a homomorphism
sending $\Sigma^{q-1}(\iota \cup Sq^2\iota)$ to the corresponding
element. It follows that $\partial = 0$

When $n= 3(mod4)$,there is a map  $g:S^{r+1} \to
\Sigma^{q-1}K(\Z,n)$ such that $g^*(\Sigma^{q-1}(\iota \cup
Sq^2\iota))\neq 0$ since otherwise $\pi_{2n}^s(K(\Z,n-1))\neq 0$.
By the above fact on map $\rho_{2^i}$ it is easy to see that there
is a map $h:S^{r+1} \to \Sigma^{q-1}K(\Z_{2^i},n)$ such that
$h^*(\Sigma^{q-1}(\iota \cup Sq^2\iota))\neq 0$.It follows from
the lemma above that $\partial h\neq 0$.
\end{proof}

With the help of Proposition\ref{T:0.2} and the known results
about $\pi _{2n+j}^s(K(\pi ,n))$ for $j=0,1$, we can now determine
the group $\pi _{2n}^s(K(\Z_{2^i} ,n-1))$.

Assume $i \geq 2$ in the following unless otherwise stated.
\begin{prop}
If $n= 0 (mod2)$,then  $${\rho_{2^i}}_*: \pi _{2n}^s(K(\Z ,n-1))
\to \pi _{2n}^s(K(\Z_{2^i} ,n-1))$$ is an isomorphism.
\end{prop}
Before the proof of the Proposition, let's give two remarks which
are clear from the proof of the Proposition.
\begin{rem}\label{T:remark}
If $i=1$, ${\rho_{2^i}}_*$ is onto.
\end{rem}
\begin{rem}\label{T:remark1}
If $n=0 (mod4)$, then the spherical cohomology class in \newline $
\pi _{2n}^s(K(\Z ,n-1))$ does not belongs to the image of the
natural map: \newline $\pi_{r}(\Sigma^qF_{n-1}(\Z)) \to \pi
_{2n}^s(K(\Z ,n-1))$.
\end{rem}
\begin{proof}
Note that we have a commutative diagram
\[
\begin{CD}
\pi_{r}(\Sigma^qF_{n-1}(\Z)) @>>> \pi _{2n}^s(K(Z ,n-1)) @>>> \pi
_{2n+1}^s(K(\Z ,n))
\\ @V{\rho_{2^i}}_*VV @V{\rho_{2^i}}_*VV   @V{\rho_{2^i}}_*VV \\
\pi_{r}(\Sigma^qF_{n-1}(\Z_{2^i})) @>>> \pi_{2n}^s(K(\Z_{2^i} ,n-1))
@>>> \pi_{2n+1}^s(K(\Z_{2^i} ,n))
\end{CD}
\]
In the above diagram,the two left horizontal maps are injective by
Lemma\ref{T:2.2}, the left vertical map is obviously an
isomorphism while the fact that the right vertical one is also an
isomorphism follows by comparing the Whitehead exact sequences of
$\Gamma_{q-1}(K(\Z,n))$ and $\Gamma_{q-1}(K(\Z_{2^i},n))$.On the
other hand ,the fact that the right horizontal map on the bottom
line is onto follows from the long exact sequence and the known
results about $\pi _{2n+j}^s(K(\Z_{2^i} ,n))$ for $j=0,1$.
\end{proof}

\begin{prop}\label{T:reduction}
For $n= 1(mod2)$, $\pi _{2n}^s(K(\Z_{2^i} ,n-1))=\pi _{2n}^s(K(\Z
,n-1))\bigoplus \Z_2$.
\end{prop}
\begin{proof}
The relevant commutative diagram in this case is
\[
\begin{CD}
\pi_{r}(\Sigma^qF_{n-1}(\Z)) @>>> \pi _{2n}^s(K(\Z ,n-1))@>>> \pi
_{2n+1}^s(K(\Z ,n))=0 \\ @V{\rho_{2^i}}_*VV @V{\rho_{2^i}}_*VV
@V{\rho_{2^i}}_*VV \\ \pi_{r}(\Sigma^qF_{n-1}(\Z_{2^i})) @>>>
\pi_{2n}^s(K(\Z_{2^i} ,n-1)) @>>> \pi_{2n+1}^s(K(\Z_{2^i} ,n))
@>\partial>>
\end{CD}
\]
By the same argument as in the last Proposition , we know the map
$\partial$ is onto. If $n = 3(mod4)$,the two left horizontal maps
are trivial by Lemma\ref{T:2.2}, thus  $\pi_{2n}^s(K(\Z_{2^i}
,n-1))\cong \text{coker}\partial\cong \Z_2 $.

If $n = 1(mod4)$, what we can get is an exact sequence
\[
0 \rightarrow \pi_{2n}^s(K(\Z ,n-1))
\rightarrow\pi_{2n}^s(K(\Z_{2^i} ,n-1)) \rightarrow \Z_2 \rightarrow
0
\]
To complete the proof, it suffices to prove the last map in the
above sequence has a section.

To do this we need another diagram
\[
\begin{CD}
 \pi _{2n}^s(K_{n-1})@>>> \pi _{2n+1}^s(K_n)  @>\partial>>
\pi_{r-1}(\Sigma^qF_{n-1}(\Z_2)) \\ @Vj_*VV  @Vj_*VV   @Vj_*VV
\\ \pi_{2n}^s(K(\Z_{2^i} ,n-1)) @>>> \pi_{2n+1}^s(K(\Z_{2^i}
,n)) @>\partial>> \pi_{r-1}(\Sigma^qF_{n-1}(\Z_{2^i}))
\end{CD}
\]
where $j:\Z_2 \to \Z_{2^i}$ is the natural inclusion.

The same argument as above combined with the proof of Theorem 10.9
in\cite{milg} shows that the two $\partial$'s are onto and $j_*$
induces an isomorphism between kernels of two $\partial$'s.
Finally we get the following diagram which gives the desired
section.
\[
\begin{CD}
 \pi _{2n}^s(K_{n-1})@>\cong>> \Z_2 \\  @Vj_*VV   @V{\cong}VV
\\ \pi_{2n}^s(K(\Z_{2^i} ,n-1)) @>>> \Z_2
\end{CD}
\]
\end{proof}

\begin{lem}\label{T:3.7}
{\it If $n$ is odd and $Sq^1_i\in H^n( K(\Z_{2^i}, n-1),
\Z_2)$ is a generator. Then $$(Sq^1_i )_*: \pi
_{2n}^s(K(\Z_{2^i}, n-1))\to \pi _{2n}^s (K(\Z_2, n))\cong \Z_2$$ is
an epimorphism.}
\end{lem}
\begin{proof} It suffices to prove that the following map
$Sq^1:K(\Z_2,n-1) \to K(\Z_2,n)$ induces an isomorphism on $2n$-th
stable homotopy group. By the calculation in Milgram's
book\cite{milg}, the first group is generated by the class
corresponding to $Sq^1(t)\bigotimes Sq^1(t)$ and the second by $s
\bigotimes s$ where $s,t$ are the fundamental classes of the
corresponding groups. Now what we want follows from the fact that
$Sq^1$ induces a homomorphism mapping $s \bigotimes s$  to
$Sq^1(t)\bigotimes Sq^1(t)$.
\end{proof}

\begin{prop}\label{T:post}
Let $\widetilde{E}_{n+q}$ ($q$ large) be the following $2$-stage
Postnikov tower. Then there is a map $f:\Sigma^qK(\Z_4,n-1)\to
\widetilde{E}_{n+q}$ such that the composite $(\Sigma^q F_{n-1}(\Z)
\to)$ $\Sigma^q K(\Z,n-1)\overset{\Sigma^q \rho_4}{\to} \Sigma^q
K(\Z_4,n-1)\to \widetilde{E}_{n+q}$ if $n=1,2(mod4)$$(\text{or, }
n=0(mod4))$
 induces an isomorphism on $\pi_r$ where $r=q+2n$ as
above.

 (1).  $n=2(mod4)$ $$
\begin{array}{ccccc}
K_{r} & \stackrel{i_2}{\longrightarrow } &\widetilde{E}_{n+q} &  &
\\ &  & \downarrow  {\Pi _2} &  &  \\ K_{r-2}\times
K_{r}& \stackrel{i_1}{\longrightarrow }  & E_{n+q} & \stackrel
{\omega _2 }{\longrightarrow} & K_{r+1}   \\ &  & \downarrow { \Pi
_1} &  &  \\ \Sigma ^{q}K(\Z_4, n-1)& \stackrel {\Sigma ^ql_{n-1} }
{\longrightarrow}& K(\Z_4,q+n-1) & \stackrel {Sq^n\times Sq^{n+2} }
{\longrightarrow} & K_{r-1}\times K_{r+1}
\end{array}
$$ where $i_1^*(\omega _2)=Sq^2Sq^1l_{r-2} +Sq^1l_{r}$.

 (2). $n=0(mod4)$
$$
\begin{array}{ccccc}
K_{r}  \stackrel{i_2}{\longrightarrow } & \widetilde{E}_{n+q} & &
\\ & \downarrow  { \Pi _2} &  &  \\ K_{r-2}
\stackrel{i_1}{\longrightarrow } & E_{n+q} & \stackrel {\omega _2
}{\longrightarrow} & K_{r+1}     \\ &   \downarrow { \Pi _1} &  &
\\ \Sigma ^qK(\Z_4, n-1)\stackrel {\Sigma
^ql_{n-1}}{\longrightarrow} & K(\Z_4,q+n-1) & \stackrel {Sq^n }
{\longrightarrow} & K_{r-1}
\end{array}
$$ where $i_1^*(\omega _2)=Sq^2Sq^1l_{r-2}$.

 (3). $n=1(mod4)$
$$
\begin{array}{ccccc}
K_{r}  \stackrel{i_2}{\longrightarrow }& \widetilde{E}_{n+q} & &
\\ & \downarrow  {\Pi _2}  & &   \\  K_{r}
\stackrel{i_1}{\longrightarrow }& E_{n+q} & \stackrel {\omega _2
}{\longrightarrow}& K_{r+1}       \\ & \downarrow  { \Pi _1}  & &
\\ \Sigma ^qK(\Z_4, n-1) \stackrel {\Sigma
^ql_{n-1}}{\longrightarrow}& K(\Z_4,q+n-1)& \stackrel { Sq^{n+1}}
{\longrightarrow} &  K_{r}
\end{array} $$
where $i_1^*(\omega _2)=Sq^2l_{r-1} $ .
\end{prop}
\begin{proof}
Denote the tower in the  Proposition by
$\widetilde{E}_{n+q}(\Z_4)$. Denote  by $\widetilde{E}_{n+q}(\Z)$ a
similar tower in which $K(\Z_4,n+q-1)$ is replaced by $K(\Z,n+q-1)$.
By Remark \ref{T:remark1}, it is easy to see that there is a map
from $\Sigma^q F_{n-1}(\Z)$ to $\widetilde{E}_{n+q}(\Z)$ which
induces an isomorphism on $\pi_r$ when $n=0(mod4)$. On the other
hand it is not difficult to see that there is a map from the tower
$\widetilde{E}_{n+q}(\Z)$ to the tower $\widetilde{E}_{n+q}(\Z_4)$
which induces an isomorphism on $\pi_r$. It remains to prove that
the natural map $\Sigma^q\iota_{n-1}:\Sigma^qK(\Z_4,n-1)\to
K(\Z_4,n+q-1)$ can be lifted to $\widetilde{E}_{n+q}(\Z_4)$ and the
lifting is compatible to that of the map
$\Sigma^q\iota_{n-1}:\Sigma^qK(\Z,n-1)\to K(\Z,n+q-1)$ to
$\widetilde{E}_{n+q}(\Z)$.

We will give a proof only for $n= 2(mod4)$ , the other cases are
similar. Consider the fiber inclusion map $h: \Sigma ^q\Gamma
_{q}\to \Sigma ^qK(Z_4, n-1)$, we have the following
Peterson-Stein formula
 $$Sq^2Sq^1Sq^n_{h} (\Sigma
^ql_{n-1})+Sq^1Sq^{n+2}_{h} (\Sigma ^ql_{n-1})=h^*\Psi(\Sigma
^ql_{n-1})\in H^{r+1}(\Sigma ^q \Gamma _{q}, \Z_2)/Q$$ where $Q=
Sq^2Sq^1(Im h^*)+Sq^1(Im h^*)=Sq^1(Im h^*)$.

By Theorem 4.6 \cite{milg} and a familiar diagram chase argument
as in the proof of Proposition 2 in Chap.16 \cite{mosh}(see also
\cite{thomas}) , we have $\Sigma ^q(\theta \otimes \theta )\in
Sq^n_h(\Sigma ^ql_{n-1})$ and $\Sigma ^q{e^2\cup(\theta \otimes
\theta )}\in Sq^{n+2}_h(\Sigma ^ql_{n-1})$. It follows easily that
$h^*\Psi(\Sigma ^ql_{n-1})=0 \in H^{r+1}(\Sigma ^q \Gamma _{q},
\Z_2)/Q$. It is not difficult to see from this and a simple
computation that $\Psi(\Sigma ^ql_{n-1})=0$ and a lifting can be
chosen such that $\omega_2$ lies in its kernel.

To complete the proof, note that , as mentioned before , there is
a commutative diagram up to homotopy
\[
\begin{CD}
\widetilde{E}_{n+q}(\Z) @>\rho_4>> \widetilde{E}_{n+q}(\Z_4) \\
 @VVV  @VVV \\
E_{n+q}(\Z) @>\rho_4>> E_{n+q}(\Z_4) \\
  @VVV  @VVV \\
  K(\Z,n+q-1) @>\rho_4>>   K(\Z_4,n+q-1)\\
  @A\Sigma ^ql_{n-1}AA  @A\Sigma ^ql_{n-1}AA  \\
\Sigma^qK(\Z,n-1) @>\rho_4>> \Sigma^qK(\Z_4,n-1)
\end{CD}
\]
The lifting from $\Sigma^qK(\Z,n-1)$ of $\Sigma^ql_{n-1}$ and the
lifting from $\Sigma^qK(\Z_4,n-1)$ of $\Sigma^ql_{n-1}$ can be made
compatible by a modification of the lifting from
$\Sigma^qK(\Z_4,n-1)$ of $\Sigma^ql_{n-1}$. The same way the
liftings  to $\widetilde{E}_{n+q}$ can also be made compatible.
Thus we have the following commutative diagram up to homotopy
which completes the proof.
\[
\begin{CD}
\widetilde{E}_{n+q}(\Z) @>\rho_4>> \widetilde{E}_{n+q}(\Z_4) \\
 @AAA  @AAA  \\
\Sigma^qK(\Z,n-1) @>\rho_4>> \Sigma^qK(\Z_4,n-1)
\end{CD}
\]
\end{proof}

\begin{rem}

 The $2^{nd}$ k-invariant $\omega _2$ in the Postnikov tower above
gives an unique secondary cohomology operator $\Psi$ (with
$Z_4$-coefficients) associated with the Adem relation
$$\begin{array}{ll} Sq^{2}Sq^{1}Sq^{n}+Sq^{1}Sq^{n+2}=0 &n=2(mod
4)\\ Sq^{2}Sq^{1}Sq^{n}=0 &n=0(mod 4)\\ Sq^{2}Sq^{n+1}=0 &n=1(mod
4)\\
\end{array}$$
Note that $E_{n+q}$ is the universal example of the operator
$\Psi$.

 By Peterson-Stein\cite{stein1}, there are operators $\Phi $ which are S-dual to
$\Psi $(which is uniquely determined by $\Psi$) so it is a
secondary operator associated with the Adem relations: $$
\begin{array}{ll}
\chi(Sq^n)Sq^3+\chi(Sq^{n+2})Sq^1+Sq^1\chi(Sq^{n+2})=0  &
n=2(\bmod 4) \\ \chi(Sq^n)Sq^3+Sq^1\chi(Sq^{n+2})=0  & n=0(\bmod
4)
\\ \chi(Sq^{n+1})Sq^2+Sq^1\chi(Sq^{n+2})=0  & n=1(\bmod 4)
\end{array}
$$ as we stated in $\S \ref{S:intro}$.
\end{rem}

\section{Proofs of Theorems 2.3, 2.6 and 2.7}\label{S:3}

\begin{proof}[Proof of Theorem \ref{T:0.3}]  First note that
 it suffices to show this for the universal
spectrum $\widetilde {W}(n)$ since the map $i: S^0\to \widetilde
{W}(n)$ factors through $i: S^0\to Y$. Notice that $H_i(\widetilde
{W}_k(n)/S^k)=0$ for $i\leq k+2$. Thus in the following proof, we
may assume that $Y_k/S^k$ satisfies the same for $k$ large.
Assuming $k$ large, without loss of generality we can assume that
$Y_k$ is a finite complex. Write  $Y^{\ast}_k$ for the $m$
$S$-dual of $Y_k$ and $g: Y^{\ast}_k \to S^{m-k}$ for the $S$-dual
of the inclusion $i: S^k \to Y_k$. Note that $g^*(\varsigma
_{S^{m-k}})\neq 0$, where $ \varsigma _{S^{m-k}}  $ is the
cohomology fundamental class of the sphere. By the $S$-duality we
get a  commutative diagram
\[
\begin{array}{ccc}
\{ S^{2n+k}, S^{k} \wedge K( \Z_4, n-1) \} &
\stackrel{i_{\ast}}{\longrightarrow} &  \{ S^{2n+k}, Y_{k} \wedge
K( \Z_4, n-1) \}   \\ \downarrow {\cong}&  & \downarrow {\cong} \\
\{ S^{2n+m}, S^{m} \wedge K( \Z_4, n-1) \} & \stackrel{g^{\ast}}{
\longrightarrow} &  \{ S^{2n+k}\wedge  Y^{\ast}_{k} , S^{m}\wedge
K(\Z_4, n-1) \}   \\
\downarrow &  & \downarrow   {q_2}_* \\ \lbrack S^{2n+m},
\widetilde{E}_{n+m} \rbrack  & \stackrel{g^{\ast}}
{\longrightarrow} & \lbrack S^{2n+k}\wedge  Y^{\ast}_{k},
\widetilde{E}_{n+m}\rbrack
\end{array}
\]
where $\widetilde {E}_{m+n}$ is the  tower in
Proposition\ref{T:post} and $q_2: S^m\land K(\Z_4, n-1)\to
\widetilde {E}_{n+m}$ is a lifting of
 $\Sigma ^ml_{n-1}$. From the
diagram above and Proposition \ref{T:post},it suffices to show
that the homomorphism $g^*$ at the bottom line is injective.From
now on we will restrict to the case $n\equiv 2(mod4)$. The other
cases are similar. Let $i_{0}: F \to \widetilde{E}_{n+m}$ be the
fibre of the composite $\Pi_1\circ \Pi_2$. Note that $F$ can be
viewed as a fibration over $K_{2n+m-2}$ with fibre $K( \Z_4, 2n+m)$
and $k$-invariant $j_{\ast}(Sq^{2}Sq^{1})(l)$; where $$j_{\ast} :
H^{m+2n+1}(-, \Z_2) \to H^{m+2n+1}(-, \Z_4)$$ is the homomorphism
induced by the inclusion $ \Z_2 \subset \Z_4$ and $l$ is the basic
class of $K_{m+2n-2}$.

 Consider
the following commutative diagrams
\[
\begin{array}{ccccc}
 & & [S^{2n+m}, F] & \stackrel{{\cong }_{\ast}}{\longrightarrow} &
[S^{2n+m}, \widetilde{E}_{n+m}]   \\
 & & \downarrow  {J:= g^*} &  & \downarrow
  {g^*} \\
\lbrack S^{2n+k} \wedge Y_{k}^{\ast}, K( \Z_4, n+m-2)\rbrack &
\stackrel{{i_{1}}_ {\ast}}{\longrightarrow} &  \lbrack S^{2n+k}
\wedge Y_{k}^{\ast},  F \rbrack &
\stackrel{{i_{0}}_{\ast}}{\longrightarrow} &  \lbrack S^{2n+k}
\wedge Y_{k}^ {\ast},  \widetilde E_{n+m}\rbrack
\end{array}
\]
and
\[
\begin{array}{ccccc}
& [S^{2n+m},  K( \Z_4, 2n+m)] & \stackrel{{\cong}}{\longrightarrow}
& [S^{2n+m},  F] & \\ & \downarrow
    {g^*}  & & \downarrow  {J}  &\\
\lbrack S^{2n+k}\land Y_k^*, K_{2n+m-3}\rbrack &
\stackrel{{j_*(Sq^2Sq^1)} }{\longrightarrow} \lbrack S^{2n+k}
\wedge Y_{k}^{\ast}, K( \Z_4, 2n+m)\rbrack & \stackrel{{\cong}
}{\longrightarrow} &  \lbrack S^{2n+k} \wedge Y_{k}^{\ast},
F\rbrack
\end{array}
\]
where $i_{1} : K(\Z_4, n+m-2) \to  F$ is the  homotopy fibre of
$i_{0}$.  $j_*(Sq^2Sq^1)$  in the second diagram above is zero
since  $Sq^3U_k=0$ and thus by duality  $\chi
(Sq^3)H^{m-k-3}(Y_k^*)=Sq^2Sq^1 H^{m-k-3}(Y_k^*)=0$. Thus the
second diagram implies that $J$ is a monomorphism. To complete the
proof, it suffices to  show $Ker(i_{0})_{\ast} = Im(i_1)_{\ast}
=0$ in the first diagram above.

 Let $q=m-n-k-1$, if $x\in H^{q-1}(Y^{\ast}_k,  \Z_4)$,then
   $Sq^n(x)\in H^{n+q-1}(Y_k^*, \Z_2) \cong
(H^{k+2}(Y_k, \Z_2))^*=0$.  On the other hand, by duality $\chi
(Sq^{n+2})U_k=0$ implies that $Sq^{n+2}H^{q-1}(Y^{\ast}_k,
\Z_2)=0$. Thus $$x \in \text{Ker}{Sq^n}\cap \text{Ker}{Sq^{n+2}}$$
Since ${ Y_k}$ is $\Phi$-orientable, i.e, $0\in\Phi (U_k)$.  By
\cite{stein1} that $0\in \Psi (x)$. Thus $x$ can be lifted to
$\widetilde{E}_{q-1} $ and so $(i_1)_{\ast}(x)=0$. This completes
the proof.
\end{proof}

For simplicity,denote  by $F_{n-1}(\Z_4)$ the space $K( \Z_4,
n-1)\wedge K( \Z_4, n-1)$ as before in the following proof.

\begin{proof}[Proof of Theorem 2.6] For $x\in H^{n-1}(M, \Z _4)$, let $f(x)=
(w\wedge id)A_\alpha (x)\in H_{2n}(K(\Z _4, n-1); \widetilde {W}(n))$. For $k$ large,
$f(x+y)$ is the following composition of maps

$S^1\land S^{2n+k}  \stackrel{id\land \Delta
\alpha}{\longrightarrow} S^1\land T\xi\wedge M_{+}
\stackrel{id\land {w\land (x\times y)}}{\longrightarrow}\newline
\rightarrow S^1\land\widetilde {W}(n)_k \wedge (K( \Z_4, n-1)\times  K( \Z_4, n-1))=
\newline = \widetilde {W}(n)_{k}\wedge S^1\land (K( \Z_4,
n-1)\times K( \Z_4, n-1)) \stackrel{id\land  \kappa
}{\longrightarrow} \widetilde {W}(n)_{k}\wedge S^1\land  K( \Z_4,
n-1), $

  \noindent where $\kappa^*(l)=l\otimes 1+1\otimes l $ for the
basic class $l\in H^{n-1}(K(\Z_4, n-1), \Z_4)$.

Identifying
$\widetilde {W}(n)_{k}\wedge S^1\land (K( \Z_4, n-1)\times K( \Z_4,
n-1))$ with $$\{\widetilde {W}(n)_{k}\wedge S^1\land K( \Z_4,
n-1)\} \lor \{ \widetilde {W}(n)_{k}\wedge S^1\land K( \Z_4,
n-1)\} \lor  $$ $$\lor \{ \widetilde {W}(n)_{k} \wedge S^1\wedge
F_{n-1}(\Z_4)\}.$$ It is readily to see that
$f(x+y)=f(x)+f(y)+g$, here $g$ is the composition

 {\small
$S^{2n+k+1} \stackrel{id \wedge \Delta \alpha}{\longrightarrow}
S^1\wedge T\xi\wedge M_{+} \stackrel{id\wedge w\land
\Delta} {\longrightarrow} S^1\land \widetilde {W}(n)_{k} \wedge
M_{+} \wedge M_{+} \stackrel{id\land x\land y}{\longrightarrow}
\widetilde {W}(n)_k \wedge S^1\land K( \Z_4, n-1) \wedge K( \Z_4,
n-1) \stackrel{id\land H(\kappa )}{\longrightarrow} \widetilde
{W}(n)_{k}\wedge S^1\land K( \Z_4, n-1), $ }
\newline where $H(\kappa )$ is the Hopf constuction of $\kappa$.

 Now the cofibration $$ S^{k+1}\wedge F_{n-1}(\Z_4)\overset{\Sigma i\wedge id}{\to} S^1\wedge \widetilde
{W}(n)_{k}\wedge F_{n-1}(\Z_4) \to S^1\wedge(\widetilde
{W}(n)_{k}/S^{k})\wedge F_{n-1}(\Z_4) $$ is also a fibration at
least in the stable range. It follows immediately that $$(\Sigma
i\wedge id)_{\ast}: \pi_{2n+k+1}(S^{k+1}\wedge F_{n-1}(\Z_4))\to
\pi _{2n+k+1}(S^1\wedge \widetilde {W}(n)_{k}\wedge F_{n-1}(\Z_4))
$$ is surjective. On the other hand, it is easy to know that the
generator $\beta \in \pi ^s_{2n}(F_{n-1}(\Z_4)) \cong \Z_2 $
satisfies $\beta ^*(l\otimes Sq^2l)\neq 0$. Thus, for the
inclusion map $i$, the composition $(\Sigma i\land id) \circ \beta
\in \pi_{2n+k+1}(S^1\wedge \widetilde {W}(n)_{k}\wedge
F_{n-1}(\Z_4))$ induces a nontrivial homomorphism on the
$(2n+k)$-th homology and thus ($\Sigma i\wedge id)_{\ast} $ is an
isomorphism. Moreover, the generator $g_0 \in \pi
^s_{2n}(\widetilde {W}(n)_k \land F_{n-1}(\Z_4)) $ satisfies that
$g_{0} ^*(U_{k}\wedge Sq^2l_{n-1} \wedge l_{n-1})\ne 0$. Thus the
composition $(id\wedge x\wedge y)( w\wedge \Delta
)(\Delta \alpha )$ is null homotopy if and only if $ \langle x\cup
Sq^{2}y, [M]_2\rangle = 0 $. By Proposition \ref{T:0.2}, the proof
now follows by the commutative diagram
\newline
\[
\begin{array}{ccc}
S^{k}\wedge \Sigma F_{n-1}(\Z_4) & \stackrel{i\wedge
id}{\longrightarrow} & \widetilde {W}(n)_{k}\wedge \Sigma
F_{n-1}(\Z_4)\\ \downarrow {id\wedge H(\kappa )} & & \downarrow {id
\wedge H(\kappa )}\\ S^{k} \wedge \Sigma K( \Z_4, n-1) &
\stackrel{i\wedge id} {\longrightarrow} & \widetilde {W}(n)_{k}
\wedge \Sigma K( \Z_4, n-1).
\end{array}
\]
\end{proof}

\begin{proof}[Proof of Theorem 2.7] Let $\mu : K_{n+k+1}\times \widetilde W_k(n)\to \widetilde
W_k(n)$ denote the fiber multiplication. Since
$d_1(w_1, w_2)=0$, $w_2$ is the composition $$T\xi \stackrel
{\Delta }{\longrightarrow} T\xi \times T\xi \stackrel {w_1\times
vU_k}{\longrightarrow} \widetilde W_k(n)\times K_{n+k+1}\stackrel
{\mu}{\longrightarrow} \widetilde W_k(n),$$ where $vU_k\in
H^{k+n+1}(T\xi , Z_2)$  is the second difference of $w_1$ and $w_2$, i.e, the secondary obstruction
to deform $w_1$ to $w_2$.  Consider the commutative diagram:
\begin{displaymath}
\begin{array}{cccc}
S^{2n+k} & \stackrel {\alpha '}{\longrightarrow} & (T\xi \land
M_+)\lor (T\xi \land M_+) & \stackrel {a}{\longrightarrow}
\widetilde {W}_k(n)\land K(Z_4, n-1)\\
\parallel &  &  {\bigcap} & \parallel \\
S^{2n+k} & \stackrel {\Delta \alpha }{\longrightarrow} & (T\xi
\times  T\xi ) \land M_+ &  \stackrel {b}{\longrightarrow}
\widetilde {W}_k(n)\land K(Z_4, n-1)
\end{array}
\end{displaymath}
where $\alpha '$ is a lifting of $\Delta \alpha$, $b=\mu (w_1\times
vU_k)\land x$, $a=(w_1\land x )\lor c$, and $c=i(vU_k)\land x$, $i:
K_{n+k+1}\to \widetilde W_k(n)$  the inclusion of the fibre.  Write
$\alpha '= \alpha _1+\alpha _2$, here $\alpha  _1$
and $\alpha _2$ are the factors of the wedge. Note that
$\phi _ 2(x)=h(b\circ \Delta \alpha )=h(a\alpha _1)+h(a\alpha _2)=
\phi _ 1(x)+h(a \alpha _2)$.  \newline
We are going to show $h(a\alpha _2)=0$.

 As $a\alpha _2$ factors through the map $i\land id: K_{n+k+1}
\land K(\Z_4, n-1) \to \widetilde W_k(n)\land K(\Z_4, n-1)$, it
suffices to prove that $$(i\land id)_*: \pi _{2n+k}(K_{n+k+1}\land
K(\Z_4, n-1)) \to \pi _{2n+k}( \widetilde W_k(n)\land K(\Z_4,
n-1))$$ is zero. Note the homomorphism $$(Sq^1\land id)_* : \pi
_{2n+k}(K_{n+k} \land K(\Z_4, n-1))\to \pi _{2n+k}(K_{n+k+1} \land
K(\Z_4, n-1)) \cong \Z_2 $$ is an isomorphism as it induces  an
ismomorphism on the $(2n+k)$-th homology groups.  The composition
$K_{n+k}\stackrel {Sq^1} {\longrightarrow} K_{n+k+1}\stackrel {i
}{\longrightarrow} \widetilde W_k(n)$ is null homotopy. Thus
$(i\land id)_*=0$. This completes the proof.
\end{proof}

\bigskip

\section{Proof of Theorem 1.9}

\vskip 4mm

In this section we prove Theorem 1.9. We first study the properties of the invariants $\mu _M$ and
$q_M(Sq^1_i)$ defined in $\S$1.

\begin{lem}\label{T:1.6}
{\it  Let $M$ be a framed manifold of dimension $2n$ with $n$ odd. Let
$q_M: H^n(M, \Z_2)\to \Z_2$ be the Kervaire quadratic form. For
$x\in H^{n-1}(M, \Z_{2^i})$,
\newline
(i) $n=3(mod 4)$, $[M, x]$ is reduced bordant to zero iff
$q_M(Sq^1_i)x=0.$\newline (ii) $n=1(mod 4)$, $[M, x]$ is reduced
bordant to $[M', x']$ where $x'\in H^ {n-1}(M')$ iff
$q_M(Sq^1_i x)=0$.}
\end{lem}
\begin{proof}
Identify the reduced framed bordism group $\widetilde{\Omega }_{2n}^{fr}(-)$ with the stable homotopy
group $\pi _{2n}^{s}(-)$. Recall that $ \pi _{2n}^s (K(\Z_2,n)) =\Z_2$. By [4] it is easy to see that
the homomorphism
$$(Sq_i^1)_*: \pi _{2n}^{s}(K(\Z _{2^i},n-1)\to \pi _{2n}^s (K(\Z_2,n))$$
is identified with the following geometrically defined homomorphism
$$
\begin{array}{lcc}
\widetilde {\Omega }_{2n}^{fr}(K(\Z _{2^i},n-1)) & \to &  \Z_2 \\
\ [M , x] & \to &  q_{M}(Sq^1_i)x
\end{array}
$$ By Theorem 3.1 and Lemma 3.7 (i) follows since $(Sq^1_i)_*$ is
an isomorphism. To prove (ii), note that there is an exact
sequence by Proposition\ref{T:reduction} and Lemma\ref{T:3.7}$$\pi
_{2n}^s(K(\Z,n-1))\to \pi _{2n}^s(K(\Z
_{2^i},n-1))\stackrel{(Sq^1_i)_*} \longrightarrow \pi
_{2n}^s(K_n).$$ This completes the proof.
\end{proof}

Now we want to study which bilinear forms $\mu$ can be realized by $(n-2)$-connected $2n$-dimensional
$\pi$-manifolds. Note that a sphere bundle over $S^{n+1}$ with fiber $S^{n-1}$ is a $\pi$-manifold
if the characteristic map of the bundle, $\theta \in \pi _n(SO(n))$, belongs to the kernel of the
stablization homomorphism $S_*: \pi _n(SO(n))\to \pi _n(SO)$.
Recall that the homotopy groups of $\pi _n(SO(n))$ are as follows
(c.f: [11]):

\vskip 3mm

\small{
 \centerline {$\pi _n(SO(n))$, $n\geq 3, \ne 6$ }
\vskip 1mm

{\footnotesize  \begin{tabular}{|r|c|c|c|c|c|c|c|c|} \hline
 $n\geq 3, \ne 6$ &  $8s $ & $8s+1$  & $8s+2$ & $8s+3$ &$8s+4$
& $8s+5$ & $ 8s+6$ & $8s+7$
\\ \hline
 $\pi _n(SO(n))$ &  $\Z_2\oplus \Z_2\oplus \Z_2   $
& $ \Z_2\oplus  \Z_2   $ &
$\Z_4   $&
$\Z   $ & $\Z_2\oplus \Z_2$& $\Z_2$ & $\Z_4$ & $ \Z$ \\ \hline
\end{tabular} }}

\noindent  and $\pi _6(SO(6))=0$.

Let $\pi : SO(n)\to S^{n-1}$ be the canoincal $SO(n-1)$-fiberation. For a $S^{n-1}$-bundle
over $S^{n+1}$ with characteristic map $\theta \in \pi _{n}(SO(n))$, say $M_\theta$, it is easy to see
that $Sq^2: H^{n-1}(M_\theta ,\Z_2)\to H^{n+1}(M_\theta , \Z_2 )$ is an isomorphism if and only if
$\pi _*(\theta )\in \pi _n(S^{n-1})=\Z_2$ is nonzero. By duality this implies that $z\cup Sq^2z=0$ for all
$z\in H^{n-1}(M_\theta , \Z_2)$ if and only if $\pi _*(\theta )=0$. The latter is equivalent to the fact of
that the bundle has a section.

\begin{lem}\label{T:0.12}
 {\it Let $M$ be a $\pi$-manifold of dimension $2n$.
Then\newline (i) $\mu _M(x, x) =0$, $\forall x\in H^{n-1}(M, \Z_2)$
if $n=2(mod 4)$. \newline (ii) $\mu _M(x, x)=0$, $\forall x\in Im
(\rho _2:T\subset H^{n-1}(M, \Z_4)\to
 H^{n-1}(M, \Z_2))$, \newline if $n$ is odd where $T$ is the set of
 elements of order $4$.\newline
(iii) If $n=0(mod 4)$, then there is a $S^{n-1}$-bundle over
$S^{n+1}$, $M$, so that $\mu _{M}(x,x)\neq 0$,
where $x\in H^{n-1}(M, \Z_2)$ is a generator.}
\end{lem}

\begin{proof} For each $x\in H^{n-1}(M, \Z_2)$, consider the reduced
bordism class $[M, x]\in \widetilde {\Omega
}_{2n}^{fr}(K_{n-1})\cong \Z_2$. It is easy to see that $x\cup Sq^2x[M]$ is
a bordism invariant. One verifies the following map defines a homomorphism
$$
\begin{array}{lcc}
\widetilde {\Omega }_{2n}^{fr}(K_{n-1})&  \to &\Z_2 \\
\ [M , x] & \to &
x\cup Sq^2x[M]
\end{array}
$$
By Remark \ref{T:remark} the reduction homomorphism
$$\widetilde {\Omega }_{2n}^{fr}(K(\Z, {n-1}))\to
\widetilde {\Omega }_{2n}^{fr}(K_{n-1})$$
is surjective if $n$ is even.

If $n=2(mod 4)$, let  $\theta\in \pi _{n}(SO(n))$ be a generator. By the tables (I)(II) of \cite{kerv1}
it follows that $\theta $ lies in the image of the inclusion map $\pi _n(SO(n-1))\to \pi _{n}(SO(n))$. By
 the remark above this implies that the sphere bundle $M_\theta $ has a section. Therefore $z\cup Sq^2z=0$
for all $z\in H^{n-1}(M_\theta , \Z_2)$. On the other hand, one can verify that $[M_\theta , z]\in
\widetilde {\Omega}_{2n}^{fr}(K_{n-1})$ is a generator if $z\in H^{n-1}(M_\theta , \Z_2)$ is nonzero.
This proves (i).

If $n=0(mod 4)$, by  \cite{kerv1} there is an element $\beta \in kerS_* : \pi
_{n}(SO(n))\to \pi _n(SO)$ so that $\pi _*(\beta )$ is nonzero. This proves (iii).

If $n$ is odd, by Lemma \ref{T:1.6} the homomorphism
$$
\begin{array}{ccc}
q(Sq^1): \widetilde {\Omega }_{2n}^{fr}(K_{n-1}) \to \Z_2\\ \ [M,
x] \to q_M(Sq^1x)
\end{array}
$$ is an isomorphism. Thus there is a $\delta \in \Z_2$ so that
$\delta q_M(Sq^1x)=x\cup Sq^2x[M]$ for all $[M, x]$.  In
particular, if $x$ can be lifted to the $\Z_4$-coefficient class
with order $4$, $Sq^1x=0$ and so $x\cup Sq^2x=0$. This completes
the proof.
\end{proof}

Now we are ready to prove Theorem 1.9.
\vskip 2mm

\begin{proof}[Proof of Theorem \ref{T:0.13}]
By Theorem 1.6 the data of invariants are homotopy invariants of the manifolds.
Thus the homotopy and homeomorphism classification of such manifolds are the same.

 There is an
isomorphism
$$\widetilde{\Omega} _{2n}^{fr}(K(H, n-1))\cong \pi _{2n}^s(K(H, n-1)).$$
Therefore from Theorem 3.1 there is a reduced framed bordism class $[M,f]\in \Omega _{2n}^{fr}(K(H, n-1))$
corresponding to the given algebraic data  $[H,\mu , \phi ]$ (resp. $[H, \mu , \phi ,\omega ]$ ) if
$n$ is even (resp. odd). This together with Lemmas 5.1 and 5.2 implies this is an 1-1 correspondence.

Add some $S^{n-1}\times S^{n+1}$ to $M$ if necessary so that $f_*: H_{n-1}(M)\to H$ is surjective. By surgery
on $M$ we may assume further that $f_*: H_{n-1}(M)\to H$ is an isomorphism and $H_{n}(M,\Q)=0$. Therefore
the data can be realized by a $(n-2)$-connected $2n$-dimensional $\pi$-manifold,
$M$, so that $H_n(M,\Q)=0$ and $\pi (M)=[H, \mu ,\phi ]$ (resp.
$[H, \mu ,\phi ,\omega ]$.

Now it suffices to prove that the map $\pi$ is injective.

Suppose that $M_i$, $i=1,2$, are two framed smooth manifolds with
the same data (for TOP manifold, the similar argument works
identically). Note that the Kervaire invariants of $M_i$ must
vanish since $H_n(M_i,\Q)=0$. By the assumption there are maps
$f_i: M_i\to K(H, n-1)$, so that $(M_1, f_1)$ and $(M_2, f_2)$ are
reduced framed bordant, where $f_i$ induces an isomorphism on the
$(n-1)$-th homology groups. Since both $M_i$ framed cobordant to
some homtotopy spheres, there is a framed homotopy sphere,
$\Sigma$, so that $(M_1, f_1)$ and $(M_2\# \Sigma , f_2)$ are
framed bordant.  By Freedman \cite{freed} or Kreck \cite{kreck} it
follows that $M_1$ and $M_2\# \Sigma $ are diffeomorphic since
$H_n(M_i, \Q)=0$. Therefore $M_1$ and $M_2$ are almost
diffeomorphic. The same argument as above applies to show that
$M_1$ and $M_2$ are homeomorphic to each other. This completes the
proof.
\end{proof}

\vskip 4mm

{\bf Acknowledgements:} This paper is a revised version of
 \cite{ff}. This work began during first author's studys and visits
at Jilin University, Nankai Institute of Mathematics,
Universit\"at Mainz, Universit\"at Bielefeld, the
Max-Planck-Institut f\"ur Mathematik and I.H.E.S. He would like to
express his sincere thanks to all of those Institutions and
to Yifeng Sun and Xueguang Zhou for their encouragements and supports,
to Matthias Kreck for teaching him his surgery theory \cite{kreck}. The second author joins the project
at the later part of
the work mainly for clarifying the argument. Part of the work was
done during his visit to Korea University. He would like to thank
Prof.Woo Mooha and Department of Mathematics Education for the
hospitality.

\end{document}